\documentclass[12pt]{article}
\usepackage{amssymb}
\usepackage{amsmath}
\usepackage{epsfig}
\usepackage{graphicx}

\newif{\ifcomentarios}
\comentariosfalse

\newtheorem{theorem}{Theorem}
\newtheorem{lemma}[theorem]{Lemma}

\newtheorem{corollary}[theorem]{Corollary}

\renewcommand{\mathbf}{\boldsymbol}

\begin{document}

\author{D. H. U. Marchetti\thanks{%
Instituto de F\'{\i}sica, Universidade de S\~{a}o Paulo, Caixa Postal 66318,
05314-970 S\~{a}o Paulo, SP, Brasil. E-mail:\textit{\ marchett@if.usp.br} }
, \  V.{\ }Sidoravicius\thanks{%
Centrum Wiskunde \& Informatica (CWI), Science Park 123, 1098 XG Amsterdam
and Instituto de Matem\'{a}tica Pura e Aplicada (IMPA), Estrada Dona
Castorina 110, Jardim Bot\^{a}nico, 22460-320 Rio de Janeiro, RJ, Brasil.
E-mail:\textit{\ vladas@impa.br}} \  and \  M. E. Vares\thanks{%
Centro Brasileiro de Pesquisas F\'{\i}sicas (CBPF), Rua Dr. Xavier Sigaud
150, 22290-180 Rio de Janeiro, RJ, Brasil. E-mail: \textit{eulalia@cbpf.br}} 
}
\title{Oriented Percolation in One--Dimensional $1 /|x-y|^{2}$\\
Percolation Models}
\date{}
\maketitle

\begin{abstract}
We consider independent edge percolation models on $\mathbb{Z}$, with edge
occupation probabilities 
\begin{equation*}
p_{\{x,y\}}=\left\{ 
\begin{array}{ll}
p & \mathrm{if}\;\left\vert x-y\right\vert =1, \\ 
1-\exp \left\{ -\beta /\left\vert x-y\right\vert ^{2}\right\} & \mathrm{\
otherwise.}%
\end{array}
\right.
\end{equation*}
We prove that oriented percolation occurs when $\beta >1$ provided $p$ is
chosen sufficiently close to $1$, answering a question posed in \cite{NS}.
The proof is based on multi-scale analysis.
\end{abstract}

\section{Introduction \label{0}}

\setcounter{equation}{0} \setcounter{theorem}{0}

It is well known that $1/r^{2}$ gives the \textquotedblleft critical"
falloff for percolation in one-dimensional long range independent edge
percolation models. Moreover, for the one dimensional Fortuin--Kasteleyn
(FK) random cluster model with weighting factor $\kappa \geq 1$ and edge
occupation probabilities of the form $p_{\{x,y\}}=f(\left\vert
x-y\right\vert )$, with $r^{2}\,f(r)\to \beta >0$ as $r \to +\infty$, it is
known that for fixed $f(j)<1,j\geq 2$ and varying $p=f(1)$, the value $\beta
^{\ast }=1$ is critical in the sense that for $\beta \leq 1$ percolation
cannot occur unless $p=1$ (see \cite{AN}), while for $\beta >1$ there is
percolation provided $p$ is sufficiently close to one (see \cite{IN} and 
\cite{M}). Such results are important in the description of the phase
transition diagram for the one-dimensional long range Ising models studied
earlier by Fr\"{o}hlich and Spencer in \cite{FS1} and for the corresponding
Potts models (\cite{ACCN},\cite{IN},\cite{M}), as these spin systems can be
constructed by a random coloring of the clusters in the FK model with $%
\kappa =2$, or $\kappa >2$ integer, respectively. For the particular case of
independent edge percolation models ($\kappa =1$) earlier results were
obtained in \cite{NS}, where it was proven that $\beta ^{\ast }\leq 1$ in
this case, and that \textit{oriented percolation} occurs when $%
\lim\limits_{x\rightarrow \infty }x^{s}\,f(x)>0$ for some $1<s<2$ and
$p$ is sufficiently close to 1. The
question whether oriented percolation occurs in the boundary case $s=2$
remained unanswered. Theorem \ref{percolation} below gives an affirmative
answer; the result is stated for the particular example of edge
probabilities in (\ref{occup}) below, and oriented percolation is shown when 
$\beta >1$ and $p<1$ is sufficiently close to one. The proofs can be easily
adapted to include any $f(\cdot )$ satisfying $\lim\limits_{x\rightarrow
\infty }x^{2}\,f(x)>1$. In this sense $\beta ^{\ast }=1$ remains critical
also for oriented percolation.

Our main result (Theorem \ref{percolation}) deals with the independent
percolation model. On the other hand, known FKG inequalities and the above
mentioned representation yield at once an application to the long range
Ising (and Potts) models, which we state as Corollary \ref{ising}. This is
the reason for the preliminaries on this more general context of FK measures.

\noindent \textbf{Preliminaries.} Consider the infinite complete graph with
set of vertices $\mathcal{V}=\mathbb{Z}$ and set of edges $\mathbb{E}
=\{\{x,y\},\;x\neq y,\;x,y\in \mathbb{Z}\}$, and let $\Omega =\{0,1\}^{
\mathbb{E}}$. One-dimensional long-range FK random cluster models with
weighting parameter $\kappa \geq 1$ are probability measures on $\sigma
(\Omega )$, the usual product $\sigma $--algebra on $\Omega $. To define
them, let us first fix $\nu $ the Bernoulli product measure on $\Omega $,
with $\nu (\omega _{\{x,y\}}=1)=p_{\{x,y\}}$ given by

\begin{equation}
p_{\{x,y\} }=\left\{ 
\begin{array}{ll}
p & \mathrm{if}\;\left\vert x-y\right\vert =1, \\ 
1-\exp \left\{ -\dfrac{\beta }{\left\vert x-y\right\vert ^{2}}\right\} & 
\mathrm{otherwise,}%
\end{array}%
\right.  \label{occup}
\end{equation}
where $0<p<1$ and $\beta >0$ are fixed parameters.

\noindent \textbf{Notation.} We write $q_{\{x,y\} }=1-p_{\{ x,y\} }$; for $%
e=\{x,y\} $ we will write $p_{e}$ instead of $p_{\{x,y\} }$, and say that $e$
\textquotedblleft is open\textquotedblright\ if $\omega _{e}=1$. The length
of an edge $e=\{x,y\}$ is $|x-y|$.

\noindent \textit{Finite volume FK measures.} Given $I\subset \mathbb{Z}$,
consider $\mathbb{E}(I)=\{\{x,y\}\in \mathbb{E}\colon x,y\in I\}$, $\Omega
_{I}=\left\{ 0,1\right\} ^{\mathbb{E}(I)}$ and $\bar{\Omega}_{I}=\left\{
0,1\right\} ^{\mathbb{E}\setminus \mathbb{E}(I^{c})}$, where $I^{c}=\mathbb{Z%
}\backslash I$. Assume that $|I|<\infty $. The corresponding finite volume 
\textit{free FK-measure} is the probability measure $\mu _{\kappa ,I}^{f}$
on $\Omega _{I}$ given by 
\begin{equation}
\mu _{\kappa ,I}^{f}(A)=\frac{\displaystyle\int_{A}\kappa ^{\mathcal{C}%
_{I}(\omega )}\nu _{I}\left( d\omega \right) }{\displaystyle\int_{\Omega
_{I}}\kappa ^{\mathcal{C}_{I}(\omega )}\nu _{I}\left( d\omega \right) }\,,
\qquad A \subset \Omega_I,  \label{mukappa}
\end{equation}%
where $\nu _{I}$ is the restriction of $\nu $ to $\Omega _{I}$, and $%
\mathcal{C}_{I}(\omega )$ denotes the number of disjoint connected
components in the graph determined by $\omega \in \Omega _{I}$ (i.e. the
graph with vertices in $I$ whose edges coincide with those $e$ such that $%
\omega _{e}=1$). The corresponding \textit{wired FK-measure} $\mu _{\kappa
,I}^{w}$ is a probability measure on $\bar{\Omega}_{I}$, defined similarly
as in (\ref{mukappa}), replacing $\nu _{I}$ by $\bar{\nu}_{I}$, the
restriction of $\nu $ to $\bar{\Omega}_{I}$ (so that $A \in \sigma(\bar{%
\Omega}_{I}$) the usual product sigma algebra), and $\mathcal{C}_{I}(\omega )
$ by $\bar{\mathcal{C}}_{I}(\omega )$, the number of disjoint connected
components intersecting $I$ in the graph with vertices in $\mathbb{Z}$
determined by $\bar{\omega}$, the configuration which extends $\omega \in 
\bar{\Omega}_{I}$ by setting $\bar{\omega}_{e}=1$ for all $e\in \mathbb{E}
(I^{c})$. Thus we may see $\mu _{\kappa,I}^{w}$ as a measure on $\Omega$,
concentrated on the configurations for which all edges in $\mathbb{E} (I^{c})
$ are open. Analogously, we may think of $\mu_{\kappa ,I}^{f}$ as a
probability measure on $\Omega$, concentrated on the configurations $\omega$
such that $\omega_{e}=0$ for any $e\in \mathbb{E}\setminus \mathbb{E}(I)$.
Keeping this in mind we have the following well-known property.

\medskip

\noindent \textit{The infinite volume limit.} On $\Omega $ we consider the
usual partial order: $\omega \leq \omega ^{\prime }$ if $\omega _{e}\leq
\omega _{e}^{\prime }$ for each $e\in \mathbb{E}$. By the FKG inequality
(see \cite{F}, \cite{ACCN}), one has if $\kappa \geq 1$ 
\begin{equation*}
\mu _{\kappa ,I}^{f}(g)\leq \mu _{\kappa ,I^{\prime }}^{f}(g)\leq \mu
_{\kappa ,I^{\prime }}^{w}(g)\leq \mu _{\kappa ,I}^{w}(g)
\end{equation*}%
for any finite intervals $I\subset I^{\prime }\subset \mathbb{Z}$, 
and any non-decreasing continuous function $g:\Omega \rightarrow \mathbb{R}$%
. Thus, as $I\nearrow \mathbb{Z}$ the limit measures $\mu _{\kappa }^{f}$
and $\mu _{\kappa }^{w}$ exist. Moreover, $\mu _{\kappa }^{f}\leq \mu
_{\kappa }^{w}$ in FKG sense.\footnote{%
That is $\mu \leq \mu ^{\prime }$ if $\mu (g)\leq \mu ^{\prime }(g)$ for any 
$g$ continuous and increasing.} If $\kappa =1$, trivially $\mu _{\kappa
}^{f}=\mu _{\kappa }^{w}=\nu $. Since $p_{\{x,y\}}\equiv f(|x-y|)$ the
measures $\mu _{\kappa }^{f}$ and $\mu _{\kappa }^{w}$ are translation
invariant; both are ergodic.

For a more general and complete discussion on the construction of random
cluster measures, including issues in the infinite volume limit for general
external conditions, see e.g. \cite{G1} and \cite{G2} (focused mostly in
short range models). This is particularly delicate when $0<\kappa <1$.

\noindent Fix $\omega \in \Omega $. An alternating sequence of vertices and
edges $x=x_{1},e_{1},x_{2},\dots ,x_{n-1},e_{n-1},x_{n}=y,\;n\geq 1,$ is
called a path connecting $x$ to $y$, and we say that the path is open if $%
\omega _{e_{i}}\equiv \omega _{\{ x_{i},x_{i+1}\} }=1$, $1\leq i\leq n-1$.
We say that $C\subset \mathbb{Z}$ is connected if for any two distinct
vertices $x,y$ in $C$ there exists an open path $\pi $ connecting them. A
maximal connected set is called an \textit{open cluster}, and $C_{x}(\omega
) $ denotes the open cluster containing $x\in \mathbb{Z}$ (we write $%
C_{x}(\omega )=\{x\}$ if $\omega _{\{x,y\} }=0$, for all $y\in \mathbb{Z}%
\setminus \{x\}$). A path $\pi =\left( x_{1},\dots ,x_{n}\right)$ connecting 
$x$ to $y,\;x<y$, is called \textit{oriented} if $x_{1}=x<x_{2}<\dots
<x_{n-1}<x_{n}=y$, and we write $x\rightsquigarrow y$ when there is an open
oriented path connecting $x$ to $y$. Analogously we define $%
C_{x}^{+}=\{y\colon x\rightsquigarrow y\}$, and the event 
\begin{equation*}
\lbrack x\rightsquigarrow \infty ]=[|C_{x}^{+}|=\infty ].
\end{equation*}


We are ready to state our main result.

\begin{theorem}
\label{percolation} For any $\beta >1$, there exist $0<p_{0}<1$ such that,
if $p>p_{0}$, then 
\begin{equation}
\nu(0\rightsquigarrow \infty )\geq 1-\epsilon  \label{perc}
\end{equation}%
holds with $\epsilon =\epsilon (p)\searrow 0$ as $p\nearrow 1$.
\end{theorem}

\noindent {\textbf{Remarks}}. \textbf{1}. Let $\kappa >1$. The statements in
Theorem 4.1 of \cite{ACCN} imply that $\nu\leq \mu_{\kappa}^{f}$ in FKG
sense, provided the probabilities $p_{\{x,y\}}$ in $\mu_{\kappa }^{f}$ are
given by (\ref{occup}) with $\beta$ replaced by $\beta^{\prime }\geq \kappa
\beta $ and $p$, writing $p= 1- e^{- \beta }$ for $|x-y|=1$, replaced
by $p^{\prime} = 1- (1-p)^{\beta^{\prime}/ \beta}$. Since $\mu_{\kappa
  }^{f}\leq \mu _{\kappa }^{w}$, the same holds as well for $\mu _{\kappa}^{w}$.

\noindent \textbf{2}. Theorem \ref{percolation} should indeed extend exactly
to the FK random cluster model with $\kappa >1$. The authors believe that
using an algebraic implementation of the multiscale analysis developed in
the present work, one should be able to obtain this extension. Nevertheless,
for the moment we do not have a full proof (\cite{MSV}).

\noindent \textbf{Some related problems.} The type of questions treated here
has various sources of interest and we mention only a couple of them, which
have to do with our own motivations. Consider the following physical
problem: take the one-dimensional Ising model with pair interactions, the
couplings decaying as the inverse of square of the distance between
vertices, at inverse temperature $\beta >1$; this is the model studied by Fr%
\"{o}hlich and Spencer (\cite{FS1}), for which a phase transition was
established. Take now the finite box $[-L,L]$ and assume the Dobrushin
boundary conditions, i.e. all spins in $(-\infty ,-L]$ will be taken as $+1$%
, and all spins in $[L,+\infty )$ will be taken as $-1$. What can we say
about the behaviour of this model when $L\rightarrow \infty $? Is there any
sort of well defined \textit{interface}? This might require a direct
analysis in terms of the spin system, but it leads to a more general
question for the FK model, regarding the behavior of connected components of
each boundary conditioned not to touch each other. (Recall that by a random
coloring of the clusters, the FK model gives origin to a spin system which
interpolates the independent percolation model ($\kappa=1$), the Ising model
($\kappa =2$) and the $q$--states Potts model ($\kappa=q>2$, integer) at
inverse temperature $\beta $ and interaction $J_{\{x,y\}}=\beta ^{-1}\log (%
\frac{1}{1-p_{\{x,y\}}})$, the representation being possible for some (but
not all) boundary conditions. For details see \cite{FK,ACCN}). Though we
still do not fully understand this problem which remains unsolved, our
results might shed some light on it. In \cite{CMPR}, the authors obtain a
more precise description for very low temperatures, using cluster expansion
techniques.

An interesting corollary of Theorem \ref{percolation} is as follows.
Consider the Ising model (with $\pm 1$--valued spins) on $\mathbb{Z}_+$,
with interaction $J_{\{x,y\}}=|x-y|^{-2}$ if $\left\vert x-y\right\vert \geq
2$ and $J_{\{x,x+1\}}=J$ at inverse temperature $\beta $. Let $%
m_{L}^{0,+}(\beta )$ denote the average spin at the origin, with
\textquotedblleft one-sided\textquotedblright\ $(+)$ boundary conditions in $%
[L,\infty )$. By the above mentioned FK representation (see e.g. \cite%
{F,ACCN,IN}), we have

\begin{equation*}
m_{L}^{0,+}(\beta )=\mu _{2,[0,L]}^{w_{r}}(0\leftrightarrow +\infty ),
\end{equation*}
where $\mu _{2,[0,L]}^{w_{r}}$ stands for the random cluster measure on $%
\{0,1\}^{\mathbb{E}(\mathbb{Z}_+)}$ with $\kappa =2$ and all the edges $%
\{x,y\}$ with $x\geq L$ and $y\geq L$ being open (wired on the right).
Together with Remark 1 following Theorem \ref{percolation}, this yields the
following

\begin{corollary}
\label{ising} For any $\beta >2$, there exist $0<p_{0}<1$ such that, if $%
p>p_{0}$, then 
\begin{equation*}
\lim_{L\rightarrow \infty }m_{L}^{0,+}(\beta )\geq \mu _{2,\mathbb{Z}%
_+}^{f}(0\rightsquigarrow \infty )\geq \nu(0\rightsquigarrow \infty )\geq
1-\epsilon
\end{equation*}%
holds for $\epsilon =\epsilon (J)\searrow 0$ as $J\nearrow \infty $.
Consequently, there exists a phase transition when the thermodynamical limit
on $\mathbb{Z}_+$ is taken with $+$ boundary conditions on the right side.
\end{corollary}

\noindent \textbf{Remarks.}\textbf{1.} In the above corollary there is
a little change of notation with respect to the previously mentioned
FK measure: the measure $\mu _{2,\mathbb{Z}_+}^{f}$ is considered here
on $\{0,1\}^{\mathbb{E}(\mathbb{Z}_+)}$.

\noindent \textbf{2.} The content of Corollary \ref{ising} is an
immediate consequence of Theorem 3.4 in \cite{IN}, with the result
holding even for $\beta > 1$. 

It is also interesting to compare the result on oriented percolation and the
previous corollary with the somehow similar question on the multiplicity of
Gibbs states for Markov chains with infinite connections, where orientation
appears naturally through the time direction. Recently Johansson and {\"{O}}%
berg \cite{JO} showed that if $g$ is a regular specification and 
\begin{equation*}
\text{var}_{k}(g)=\sup\{\|g(\sigma)-g(\sigma^\prime)\|_1 \colon
\sigma_i=\sigma^\prime_i, i=1,\dots,k\},
\end{equation*}
then $g$ admits a unique Gibbs measure whenever the sequence $\{{\text{var}}%
_{k}(g)\}_{k=1}^{+\infty }$ is in $\ell ^{2} $. This tells, in particular,
that there are no multiple limiting measures for chains with connections
decaying as $r^{-2}$, as in Example 1 in \cite{JO}. This contrasts with the
two-sided Ising models and, as our Theorem says, with percolation models.
The understanding of Markov chains with infinite connections in the
non-uniqueness regime is still very poor, and it is known as a notoriously
difficult problem. There is strong evidence (see \cite{BHS}) that
multi-scale analysis techniques analogous to those developed in this work
could be turned into a robust tool to study this question.

\medskip

\noindent \textbf{Heuristics of the proof.} The proof relies on Fr\"{o}%
hlich-Spencer multi-scale analysis ideas (\cite{FS}, \cite{FS1}), and we use
the version developed in \cite{KMP} and \cite{M}. In the next few paragraphs
we outline the scheme of the proof, and comment on some key ideas, avoiding
most of consuming technical points. Our goal here is only to give a very
schematic and approximate picture, postponing precise formulations (which
tend to be quite involved) to later in the text.

\medskip

\noindent \textbf{The goal.} We look for an event of positive probability,
whose occurrence implies not only the existence of an infinite open
component, but also guarantees the presence of an oriented infinite open
path. Essentially, we will construct such an event, and show that it has
positive probability. Our key estimate will be: if $\beta >1$, we can find $%
\delta >0$, $\delta^\prime <1$ and $p$ sufficiently close to $1$, so that 
\begin{equation}
\nu(\exists \,\text{open path } \pi=(x_1,e_1,\dots,x_n)\colon x_1 \le
-L+L^{\delta^\prime},\, x_n \ge L-L^{\delta^\prime}, 0< x_i-x_{i-1} \le
L^{\delta^\prime},\, \forall i\,)\ge 1-2L^{-\delta}  \label{basica}
\end{equation}%
for $L=l_k$ as defined below and any $k \ge 1$, $l_1$ being sufficiently
large, and where $\nu$ stands for the product measure defined before. We
will have little control on how close to one $p $ has to be (or,
equivalently, on how large we need $l_1$). 

\medskip

\noindent \textbf{Scales.} We choose super-exponentially fast growing
scales. Given $1<\alpha <2$, $l_{0}=1$ and $l_{1}$ an integer sufficiently
large, let 
\begin{equation}
l_{k}=\lfloor l_{k-1}^{\alpha -1}\rfloor l_{k-1},\quad \;k=2,3,\ldots ,
\label{lk}
\end{equation}%
where as usual $\lfloor z\rfloor =\max \{n\in \mathbb{N}\colon n\leq z\}$.
We will use the so--called dynamical blocking argument, where the size and
location of blocks\footnote{%
Successive blocks share an end-vertex.} will be defined along the procedure
and will depend on the configuration at lower scales. Still, the length of
each block $I^{(k)}$ of the $k$--th level (called $k$--block) will be of
order $l_{k}$. More precisely, we shall see that $l_{k}-2l_{k}^{\alpha
^{\prime }/\alpha }-6l_{k-1}\leq |I^{(k)}|\leq 3l_{k}+6l_{k-1}$, for
suitable $1<\alpha ^{\prime }<\alpha $. (In particular, $|I^{(k)}|\ll
l_{k+1}^{\alpha ^{\prime }/\alpha }\ll l_{k+1}$, if $k\geq 1$ and $l_{1}$ is
large.)

\noindent \textbf{Defected and good blocks.} Further we will use the
following recursive definition of ``defected'' block. Fix $1 < \alpha^\prime
< \alpha$ to be specified later.

\textbf{1)} We say that the $0$--block $[i, i+1]$ is \textit{defected} if
the corresponding nearest neighbor edge $\{ i, i+1 \}$ is closed; otherwise
the $0$--block is said to be \textit{good} and the open nearest neighbor
path from $i$ to $i+1$ is called a $0$--\textit{pedestal};

\textbf{2)} For $k\ge 1$, a $k$--block $I^{(k)}=[s, s^\prime]$ is \textit{%
defected} if either it contains two or more defected $(k-1)$--blocks, or it
contains only one defected $(k-1)$--block $[i, i^\prime]$ but there is no
open edge $\{ a, a^\prime \}$ of length at most $l_k^{\alpha^\prime /
\alpha} $ , with $a\le i, \; i^\prime \le a^\prime, a\in \Upsilon,\;
a^\prime \in \Upsilon^\prime$, for some $(k-1)$--pedestals $\Upsilon, \;
\Upsilon^\prime$ contained in $I^{(k)}$. Otherwise $I^{(k)}$ is called 
\textit{good}.

Thus, if a $k$--block $[s, s^\prime]$ is good, then it contains an oriented
open path going from $s$ to $s^\prime$: in the case it has no defected $%
(k-1) $--blocks, this path can be obtained by concatenating $(k-1)$%
--pedestals of the good $(k-1)$--blocks which constitute the given $k$%
--block; if it has a (single) defected $(k-1)$--block, a similar
concatenation yields an oriented open path going from $s$ to $a$, which is
followed by an open edge $\{ a, a^\prime \} $, and then followed by another
concatenation of $(k-1)$--pedestals of good $(k-1)$--blocks, from $a^\prime$
to $s^\prime$. In both cases, such path from $s$ to $s^\prime $ will be
called $k$--\textit{pedestal}, and denoted by $\Upsilon$. The part of the
cluster between $a$ and $a^\prime$ is again disregarded in the future
construction since we have little control on oriented connectivity in this
segment. The condition $a^\prime -a\le l_k^{\alpha^\prime / \alpha}$ will be
crucial to guarantee that pedestals are quite dense sets (within the
corresponding good blocks), used to push the construction to higher levels.
Some care is needed when treating defects close to the boundary, which we
have disregarded here.

\noindent \textbf{Strategy.} Being ``defected'' doesn't necessarily imply
that there is no oriented open path connecting the endpoints of the block.
Nevertheless, in order to avoid substantial technical difficulties, we will
follow two rules that simplify our construction:

a) once a block is defected, we will assume the worst possible situation,
namely it will be considered as if all edges within this block were closed.

b) once we have at least two defected $(k-1)$--blocks within a $k$--block $I$%
, we will not try to find connections within the $k$--block to fix its
connectivity, but rather will ``push the problem to the next level'', and
try to ``jump over'' this troubled block $I$ by a longer edge of length at
most $l_{k+1}^{\alpha^\prime/ \alpha}$, which starts at the pedestal of some
good $k$--block to the left of $I$, and ends similarly on 
the right of $I$.

\noindent \textbf{Estimates.} The scale $l_{1}$ will be taken large enough,
to be determined later depending on the parameter $\beta >1$ and the
auxiliary parameters $\delta >0$, $1<\alpha ^{\prime }<\alpha <2$, to be
chosen at the end of Sec.2 (see (\ref{aa})--(\ref{aald})). Once $l_{1}$ is
chosen, we shall take $p$ so that: 
\begin{equation}
p\geq \left( 1+\frac{(\ln 2)^{5}}{128}\,l_{1}^{-\delta -1}\right) ^{-1}.
\label{choice00}
\end{equation}

For $k\geq 2$, let $I^{(k)}$ be a $k$--block of length\footnote{%
This is not exact in general, but holds approximately, cf. (\ref{vno}).} $%
l_{k}$, which consists of $N_{k}=l_{k}/l_{k-1}=\lfloor l_{k-1}^{\alpha
-1}\rfloor $ blocks of level $(k-1)$, of length $l_{k-1}$, and written as $%
\{I_{j}^{(k-1)}\}_{j=1}^{N_{k}}$. Assume that we have the following estimate 
\begin{equation*}
\nu (I_{j}^{(k-1)}\;{\text{is defected}})\leq l_{k-1}^{-\delta },\quad 1\leq
j\leq N_{k}.
\end{equation*}
Under the above assumptions, and if $\delta $ is chosen to satisfy (\ref{del}%
), we see that 
\begin{equation}
\nu \left( \exists \;1\leq i<j\leq N_{k}\colon I_{i}^{(k-1)},I_{j}^{(k-1)}%
\text{ are both defected}\right) \leq \frac{1}{2}\,l_{k}^{-\delta }.
\label{psiu}
\end{equation}
When the defected $I_{i}^{(k-1)}$ is unique, we assume for the moment that
it stays at distance larger than $l_{k}^{\alpha ^{\prime }/\alpha }$ from
the boundary of $I^{(k)}$. (Otherwise a sequence of local adjustments of
blocks will be needed, as we shall see in Sect. \ref{I}. The left- and
right-most extremal blocks in our volume are treated differently.) In this
case let $a$ and $a^{\prime }$ be the end-vertices of the unique defected
block $I_{i}^{(k-1)}$. By our construction, there exists an oriented path
starting from the left boundary of $I^{(k)}$ and ending at the vertex $a$
and another open oriented path starting from vertex $a^{\prime }$ and going
to the right boundary of $I^{(k)}$. Both these paths are obtained by
concatenating pedestals of all good $(k-1)$--blocks on the left side of the
defected block $I_{i}^{(k-1)}$ and, respectively, on the right side. We
denote these new left and right pedestals by $\Upsilon $ and $\Upsilon
^{\prime }$, respectively. Given that $I^{(k)}$ has a unique defected $%
I_{i}^{(k-1)}=[a,a^{\prime }]$, and given the pedestals $\Upsilon $ and $%
\Upsilon ^{\prime }$, one has the following upper bound for the conditional $%
\nu $--probability of not finding an open edge $\{x,y\}$ with $x\leq
a,a^{\prime }\leq y,\;x\in \Upsilon ,\,y\in \Upsilon ^{\prime }$ and $%
y-x\leq l_{k}^{\alpha ^{\prime }/\alpha }$ : 
\begin{equation}
\prod_{\substack{ x,y:x\leq a\,<a^{\prime }\leq y,  \\ y-x<l_{k}^{\alpha
^{\prime }/\alpha }  \\ x\in \Upsilon ,\,y\in \Upsilon ^{\prime }}}%
q_{\{x,y\}}\,=\exp \Big\{-\sum_{\substack{ x,y:x\leq a\,<a^{\prime }\leq y, 
\\ y-x<l_{k}^{\alpha ^{\prime }/\alpha }  \\ x\in \Upsilon ,\,y\in \Upsilon
^{\prime }}}\frac{\beta }{\left\vert x-y\right\vert ^{2}}\Big\}\leq
l_{k-1}^{-\beta (1-\eta )(\alpha ^{\prime }-1)},  \label{conserto}
\end{equation}%
where $\eta =\eta (\alpha ,\alpha ^{\prime },l_{1})>0$, and can be taken
arbitrarily small if $l_{1}\rightarrow \infty $.

The precise statement and proof of the above estimate will be given in Lemma %
\ref{est33}. It requires some work, and in order to obtain it for suitable $%
\eta =\eta (\alpha ,\alpha ^{\prime },l_{1})>0$ which can be taken
arbitrarily small if $l_{1}\rightarrow \infty $ we will need to use certain
geometric properties of pedestals $\Upsilon $ and $\Upsilon ^{\prime }$,
which propagate inductively from each level into the next one. Namely, the
pedestals are relatively dense sets (see (\ref{bagunca}) in Sect. \ref{I})
as the construction will show. Using the above estimate, writing 
\begin{eqnarray*}
&&\text{$\big\{I^{(k)}$ has a unique defected $(k-1)$--block $[a,a^\prime]$
and remains defected $\big\}=$} \\
&&\text{$\big\{I^{(k)}$ has unique defected $(k-1)$--block $[a,a^{\prime }]%
\big\}\cap $} \\
&&{\text{$\big\{$ there is no open edge $\{ x,y\} $ with $x\leq a,a^{\prime
}\leq y,\;x\in \Upsilon ,\,y\in \Upsilon ^{\prime }$ and $y-x\leq
l_{k}^{\alpha ^{\prime }/\alpha }$ $\big\}$ }},
\end{eqnarray*}%
and since these events depend on disjoint sets of edges, we easily get: 
\begin{equation}
\nu \text{$\big(I^{(k)}$ has a unique defected $I_{i}^{(k-1)}$ and remains
defected}\big)\leq l_{k-1}^{\alpha -1-\delta }l_{k-1}^{-\beta (1-\eta
)(\alpha ^{\prime }-1)}\leq \frac12\,l_{k}^{-\delta },  \label{aaa}
\end{equation}
provided 
\begin{equation}
\beta (1-\eta )(\alpha ^{\prime }-1)>(\delta +1)(\alpha -1).  \label{delta-2}
\end{equation}
Since $\beta >1$ and $\eta =\eta (\alpha ,\alpha ^{\prime },l_{1})$ can be
taken very small provided $l_{1}$ is large, it will suffice to suitably fix
the parameters $\alpha $ and $\alpha ^{\prime }$ ($\alpha ^{\prime }$ close
enough to $\alpha $). This is done at the end of Sect. \ref{I}.

\medskip

\noindent \textbf{Difficulties.} To carry on this scheme we have to go
through several \textquotedblleft unpleasant\textquotedblright\ and rather
involved points. The use of a dynamical blocking argument, with the blocks
of a given level depending not only on the size and location of lower level
blocks, but also on their \textquotedblleft status\textquotedblright\
(defected or good), requires a rather tight bookkeeping. This is expressed
through what we call \textquotedblleft itineraries\textquotedblright .

Once this is achieved, all necessary estimates follow along the scheme of 
\cite{FS1} and \cite{KMP}.

In the next section we define the blocks and describe the dynamic
renormalization procedure, proving Theorem \ref{percolation}.





\section{Spatial blocks (Dynamic Renormalization) \label{I}}

\setcounter{equation}{0} \setcounter{theorem}{0}

\noindent \textbf{Notation.} For $L\in {\mathbb{N}}$, assumed to be large,
the construction will involve the configuration $\omega $ restricted to the
set of edges with both end-vertices in $[-L,L]$, where $[a,b]=[a,b]\cap 
\mathbb{Z}$ throughout.\footnote{%
Except in the proof of Lemma 2.1.} We write $\Omega _{L}$ as a shorthand for 
$\Omega _{\lbrack -L,L]}$. Scales $\{l_{k}\}_{k\in {\mathbb{N}}}$ are
defined in the following way: $l_{0}=1$, given $\beta >1$ we shall take
auxiliary parameters $\delta >0$, $1<\alpha ^{\prime }<\alpha <2$ chosen
according to (\ref{aa})--(\ref{aald}), $l_{1}$ will be a suitably large
integer and the parameter $p<1$ will be taken sufficiently close to $1$,
depending on $l_{1}$. Then we let $l_{k}$ be given by (\ref{lk}).

Further we denote $x_{j}^{(k)}=jl_{k},\;j\in {\mathbb{Z}}$.

\noindent For the proof of Theorem \ref{percolation} we may assume that $L =
l_M$, for some $M \in {\mathbb{N}}$.

\medskip

\noindent Throughout the text ${\mathbb{I}}_A$ stands for the indicator
function of an event $A$, i.e. ${\mathbb{I}}_A(\omega)=1$ or 0 according to $%
\omega \in A$ or not.

\bigskip

\noindent \textbf{Decomposition of events. Level 0.} We set $I_i^{(0)} = [i,
i+1]$. They are called $0$--blocks, and for $i$ such that $I_i^{(0)} \subset
[-L,L] $ we define the events: 
\begin{equation*}
G(I_i^{(0)} )= \{ \omega : \omega_{\{i,i+1\}}=1\}, \qquad B(I_i^{(0)} )= \{
\omega : \omega_{\{i,i+1\}}=0\}.
\end{equation*}
$I_{i}^{(0)}$ is said to be defected when $B(I_{i}^{(0)})$ occurs; otherwise
it is said to be a good $0$--block. \medskip

\noindent \textbf{Level 1.} Consider the intervals $\widetilde{I}
_{j}^{(1)}=[jl_{1},(j+1)l_{1}]$ and for each $j$ such that $\widetilde{I}
_{j}^{(1)}\subset [-L,L]$ we define the following partition of $\Omega _{L}$%
: 
\begin{eqnarray}
G(\widetilde{I}_{j}^{(1)})
&=&\bigcap_{i=jl_{1}}^{(j+1)l_{1}-1}G(I_{i}^{(0)}),  \notag \\
H_{i}(\widetilde{I}_{j}^{(1)})&=& B(I_{i}^{(0)})\cap \bigcap\limits 
_{\substack{ {s=jl_{1}}  \\ {s\neq i}}}^{(j+1)l_{1}-1}G(I_{s}^{(0)})\quad {%
\text{ for $\; i \in [jl_{1}, (j+1)l_{1}-1]$}},  \notag \\
H(\widetilde{I}_{j}^{(1)}) &=&\bigcup_{i=jl_{1}}^{(j+1)l_{1}-1}H_{i}(%
\widetilde{I}_{j}^{(1)}),  \notag \\
B(\widetilde{I}_{j}^{(1)}) &=&\left( G(\widetilde{I}_{j}^{(1)})\cup H(%
\widetilde{I}_{j}^{(1)})\right) ^{c},  \label{eventos111}
\end{eqnarray}
where $G$ stands for \textit{good}, $H$ for \textit{hopeful} and $B$ for 
\textit{bad}, and accordingly, $\widetilde{I}_{j}^{(1)}$ is said to be good
(for given $\omega$) if it contains no defected $0$--blocks, ``hopeful" if
it contains only one defected $0$--block, and is said to be ``bad"
otherwise. When $H_{i}(\widetilde{I}_{j}^{(1)})$ occurs, $I_{i}^{(0)}$ is
called the defected 0-block in $\widetilde{I}_{j}^{(1)}$.

\noindent \textbf{Adjustment.} Given $\omega $, we first consider the set of
all $j$'s such that $\omega \in H_{i_j}(\widetilde{I}_{j}^{(1)})\subset H(%
\widetilde{I}_{j}^{(1)})$ and such that the index $i_j$ of the (unique)
defected block $I_{i_{j}}^{(0)}\subset \widetilde{I}_{j}^{(1)}$ verifies $%
jl_{1}\leq i_{j}\leq jl_{1}+\lfloor l_{1}^{\alpha ^{\prime }/\alpha }\rfloor
-1$ (resp.\ $(j+1)l_{1}-\lfloor l_{1}^{\alpha ^{\prime }/\alpha }\rfloor
\leq i_{j}\leq (j+1)l_{1}-1$).

If this set is empty in both cases, we set $I^{(1)}_j = \widetilde I^{(1)}_j 
$ for all $j$'s, and say that $G(I_{j}^{(1)})$, $H(I_{j}^{(1)})$, $%
B(I_{j}^{(1)})$ occurs, according to the occurrence of the corresponding $G(%
\widetilde{I}_{j}^{(1)})$, $H(\widetilde{I}_{j}^{(1)})$, $B(\widetilde{I}%
_{j}^{(1)})$.

If this set is not empty, we take arbitrarily one of such indices $j$; if $%
\widetilde I^{(1)}_j $ is not the interval which contains $-L$ (resp. $L$),
to be treated in case 3) below, we check if $\widetilde{I}_{j-1}^{(1)}$
(resp.\ $\widetilde{I}_{j+1}^{(1)}$) has a defected $0$--block in the
sub-interval $[jl_{1}-2\lfloor l_{1}^{\alpha ^{\prime }/\alpha }\rfloor
,jl_{1}]$ (resp.\ $[(j+1)l_{1},(j+1)l_{1}+2\lfloor l_{1}^{\alpha ^{\prime
}/\alpha }\rfloor -1]$).

1) If yes, then we consider a new interval $I_{j-1}^{(1)} = \widetilde
I_{j-1}^{(1)}\cup \widetilde I_j^{(1)}$ (resp.\ $I_{j}^{(1)} = \widetilde
I_{j}^{(1)}\cup \widetilde I_{j+1}^{(1)}$) and say that the event $%
B(I_{j-1}^{(1)} )$ (resp.\ $B(I_{j}^{(1)} )$) occurs. (This is motivated by
the fact that for the chosen $\omega$ the new interval will contain at least
two defected $0$--blocks.)

\medskip

2) If not, then we consider two new intervals $%
I_{j-1}^{(1)}=[(j-1)l_{1},jl_{1}-\lfloor l_{1}^{\alpha ^{\prime }/\alpha
}\rfloor ]$ and $I_{j}^{(1)}=[jl_{1}-\lfloor l_{1}^{\alpha ^{\prime }/\alpha
}\rfloor ,(j+1)l_{1}]$ (resp.\ $I_{j}^{(1)}=[jl_{1},(j+1)l_{1}+\lfloor
l_{1}^{\alpha ^{\prime }/\alpha }\rfloor ]$ and $I_{j+1}^{(1)}=[(j+1)l_{1}+%
\lfloor l_{1}^{\alpha ^{\prime }/\alpha }\rfloor ,(j+2)l_{1}]$). We say that 
$H(I_{j}^{(1)})$ occurs, and that $G(I_{j-1}^{(1)})$, $H(I_{j-1}^{(1)})$, $%
B(I_{j-1}^{(1)})$ occurs according to the occurrence of the corresponding
event $G(\widetilde{I}_{j-1}^{(1)})$, $H(\widetilde{I}_{j-1}^{(1)})$, $B(%
\widetilde{I}_{j-1}^{(1)})$ (resp. we say that $H(I_{j}^{(1)})$ occurs, and
that \ $G(I_{j+1}^{(1)})$, $H(I_{j+1}^{(1)})$, $B(I_{j+1}^{(1)})$ occurs
according to the occurrence of the corresponding event \ $G(\widetilde{I}%
_{j+1}^{(1)})$, $H(\widetilde{I}_{j+1}^{(1)})$ $B(\widetilde{I}_{j+1}^{(1)})$%
). In this case the adjustment moves the boundary \textquotedblleft
away\textquotedblright\ from the unique defected block in $I_{j}^{(1)}$, but
doesn't change the number of the defected $0$--blocks in the adjusted
intervals.

\medskip

3) If the interval $\widetilde I_{j}^{(1)}$ under consideration is the
leftmost (resp. the rightmost) interval in $[-L,L]$, and the defect stays
within distance less than $\lfloor l_{1}^{\alpha ^{\prime }/\alpha }\rfloor$
from $-L$ (resp. $L$), we still set $I_{j}^{(1)} = \widetilde I_{j}^{(1)}$
and say that $G(I_{j}^{(1)})$ occurs.

\medskip

4) We set $I^{(1)}_j = \widetilde I^{(1)}_j $ if $\widetilde I^{(1)}_j $ was
not involved in the previous adjustment, and say that $G(I_{j}^{(1)})$, $%
H(I_{j}^{(1)})$, $B(I_{j}^{(1)})$ occurs if, accordingly, \ $G(\widetilde{I}%
_{j}^{(1)})$, $H(\widetilde{I}_{j}^{(1)})$, $B(\widetilde{I}_{j}^{(1)})$
occurs.

To conclude this step, we re-numerate the intervals from left to right as $%
I_j^{(1)}\, j=1,\dots $. If we are still left with intervals $I_j^{(1)}$ for
which $H( I_j^{(1)} )$ occurs and its defected $0$--block stays within
distance $\lfloor l_1^{\alpha^{\prime}/ \alpha}\rfloor$ from the boundary of 
$I_j^{(1)} $, we repeat the above procedure to the intervals already
adjusted in the previous step. After finitely many steps of such adjustment
procedure there are left no intervals $I_j^{(1)}$ for which the event $%
H(I_j^{(1)})$ occurs and its defected $0$--block stays within distance $%
\lfloor l_1^{\alpha^{\prime}/ \alpha}\rfloor $ from the boundary, and the
adjustment procedure is then stopped. (Of course, due to item 3, the left-
or rightmost intervals can stay with a unique defect, if this is close
enough to $-L$ or $L$ respectively.)

\medskip

\noindent \textbf{Remark.} Notice that the adjustment procedure is well
defined, i.e. the final partition does not depend on the order in which we
do adjustments and in which order we pick the intervals that still need to
be adjusted (in case we have more than one). It also has a \textit{locality}
property, i.e. the final modification of each initial interval $\widetilde
I^{(1)}_j$ depends on the values of the configuration in the nearest
neighbor and, at most, in the next nearest neighbor intervals only.

\medskip

Once the adjustment is completed, the obtained intervals, always
re-numerated from left to right as $I^{(1)}_j,\, j=1,\dots$, are called $1$%
-blocks. Notice that $l_1 - 2\lfloor l_1^{\alpha^{\prime}/ \alpha}\rfloor
\leq | I^{(1)}_j | \leq 3l_1$, and $\cup_j I^{(1)}_j=[-L,L]$. \medskip

In other words, the restriction of $\omega $ to nearest neighbor edges of $%
[-L,L]$ determines through the above procedure a random ``partition" $%
I^{(1)}(\omega)\equiv \{I_{j}^{(1)}(\omega )\}_{j}$ of the interval $[-L,L]$
into $1$--blocks, with the property that any two adjacent blocks share an
end-vertex. This is the final state of the ``adjustment" procedure. Values
of $\omega $ on the nearest neighbor edges in $[-L,L]$ also determine where
the defected 0-blocks are located within each 1-block, and we denote by $%
D_{j}^{(1)}(\omega)$ the set of indices of the defected $0$--blocks within $%
I_{j}^{(1)}(\omega)$, and $D^{(1)}(\omega)\equiv \{D_{j}^{(1)}(\omega
)\}_{j} $. The random object $J_{L}^{(1)}:=\{I_{j}^{(1)},D_{j}^{(1)}\}$ is
called itinerary at level 1 or 1-itinerary. \medskip

\noindent \textbf{1-Pedestals.} Given the $1$--itinerary $J^{(1)}_L$, we
shall attribute to each random block $I^{(1)}_j$ a state $G$ or $B$. We
first consider the case that $I^{(1)}_j$ is not the leftmost (i.e. $j\neq 1$%
) nor the rightmost $1$--block, to be treated at the end. When $%
D_{j}^{(1)}=\emptyset$, so that all nearest neighbor edges are open, we say
that $I^{(1)}_j$ is in state $G$, and we define the pedestal $\Upsilon
(I_{j}^{(1)})=I_{j}^{(1)}$. When $|D_{j}^{(1)}| = 1 $, the set of vertices $%
x\in I_{j}^{(1)}$ to the left (resp. right) of the defected $0$--block in $%
I_{j}^{(1)}$ will be called \textit{left 1-pedestal of $I_{j}^{(1)}$} (resp. 
\textit{right 1-pedestal}) and denoted by $\Upsilon _{\mathcal{L}
}(I_{j}^{(1)})$ (resp. $\Upsilon _{\mathcal{R}}(I_{j}^{(1)})$). The vertices
in each of these 1-pedestals are connected by open nearest neighbor edges.
In this situation we say that $I^{(1)}_j$ is in state $G$ when the following
event occurs: 
\begin{equation}
\lbrack \omega :\exists x\in \Upsilon _{\mathcal{L}}(I_{j}^{(1)}(\omega)),y%
\in \Upsilon _{\mathcal{R}}(I_{j}^{(1)}(\omega)),\;1< y-x\leq \lfloor
l_{1}^{\alpha ^{\prime }/\alpha }\rfloor :\omega _{\{x,y\}}=1],
\label{2.8888}
\end{equation}
and otherwise we say that $I^{(1)}_j$ is in state $B$. Similarly, if $%
|D_{j}^{(1)}| >1$ the block $I^{(1)}_j$ is in state $B$.

For the leftmost (rightmost) $1$--block, there is some little difference: In
the case $|D_{j}^{(1)}| = 1 $ and if the unique defected $0$--block stays
within distance $\lfloor l_{1}^{\alpha ^{\prime }/\alpha }\rfloor$ from $-L$
($L$), the block is said to be in state $G$, and the pedestal is defined as
the previously defined right 1-pedestal (left 1-pedestal, resp.), $\Upsilon
(I_{1}^{(1)})= \Upsilon _{\mathcal{R} }(I_{1}^{(1)})$ ($\Upsilon
(I_{j}^{(1)})= \Upsilon _{\mathcal{L} }(I_{j}^{(1)})$, resp.). Except for
this, the definition goes as with the other blocks.

With a little abuse of notation we use again the symbols $G(I_{j}^{(1)})$
and $B(I_{j}^{(1)})$ to denote that $I_{j}^{(1)}$ is in state $G$ and $B$
respectively. We say that $I_{j}^{(1)}(\omega)$ is \textit{defected} if and
only if it is in state $B$.

In (\ref{2.8888}), if the pair $(x,y)$ such that $x\in \Upsilon _{\mathcal{L}%
}(I_{j}^{(1)}),y\in \Upsilon _{\mathcal{R}}(I_{j}^{(1)}),\;y-x\leq \lfloor
l_{1}^{\alpha ^{\prime }/\alpha }\rfloor$, $\omega _{\{x,y\}}=1$ is not
unique, we choose one in arbitrary way, and, once the pair $(x,y)$ is
chosen, the interval $[x+1,y-1]$ will be called defected part of $%
I_{j}^{(1)} $, and denoted by $\mathcal{D}(I_{j}^{(1)})$. In this case we
define $\Upsilon (I_{j}^{(1)})=\left( \Upsilon _{\mathcal{L}%
}(I_{j}^{(1)})\cup \Upsilon _{\mathcal{R}}(I_{j}^{(1)})\right) \setminus 
\mathcal{D}(I_{j}^{(1)})$.

In particular, a $1$-pedestal $\Upsilon (I_{j}^{(1)})$ is given by the
vertices of an open oriented path with all edges, except possibly one, being
nearest neighbor, and this larger edge has length at most $\lfloor
l_{1}^{\alpha ^{\prime }/\alpha }\rfloor $. For each $1$--block, except
possibly the two which contain the extremes $-L$ or $L$, the pedestal
connects left and right endpoints of the interval. In the leftmost
(rightmost) case, it is allowed for the $1$--pedestal to start (end) at a
vertex within distance $\lfloor l_{1}^{\alpha ^{\prime }/\alpha }\rfloor +1 $
of $-L$ ($L$ respectively).



\begin{figure}[tbp]
\centering
\includegraphics[scale=0.35]{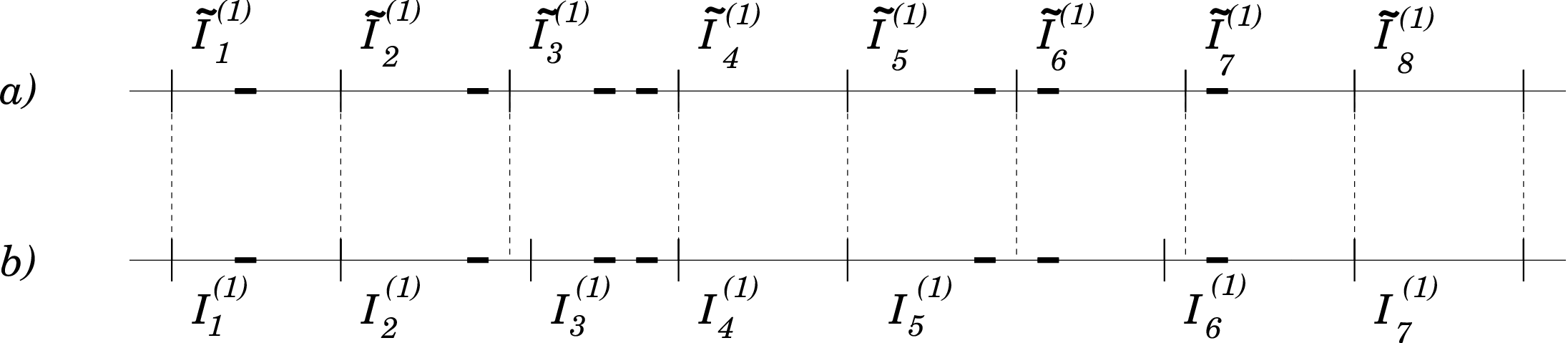}
\caption{Adjustments: part $a)$ shows the deterministic $1$--blocks $%
\widetilde{I}_{j}^{(1)}$, bold-face segments show location of the defects.
Part $b)$ shows how these blocks were adjusted. $\widetilde{I}_{5}^{(1)}$
and $\widetilde{I}_{6}^{(1)}$, merge into a single $1$--block $I_{5}^{(1)}$. 
}
\label{fig1}
\end{figure}


\noindent \textbf{Level {\textit{k} }.} Let $2\le k\le M$. Assume to have
completed the step $(k-1)$ of the recursion. In particular, for each $\omega
\in \Omega_L$ and any $r=1,\dots ,k-1$ the following objects are defined:

\begin{itemize}
\item the collection of $r$--blocks $I^{(r)}(\omega)=\{I_{j}^{(r)}(\omega)
\}_{j} $, such that $\cup_{j}I_{j}^{(r)}(\omega)=[-L,L]$, and any two
adjacent intervals share exactly an endpoint. Moreover, the uniform bound
holds: 
\begin{eqnarray}
l_{r}-(2\lfloor l_{r}^{\alpha ^{\prime }/\alpha }\rfloor
+6l_{r-1})<|I_{j}^{(r)}(\omega)|\leq 3l_{r}+6l_{r-1},  \label{comprimento}
\end{eqnarray}

\item each of the $I_{j}^{(r)}(\omega)$ can be in two possible states $G$ or 
$B$: \newline
If $I_{j}^{(r)}(\omega)$ is in state $G$ and it is not the leftmost or the
rightmost interval of the partition, then $\omega$ has an $r$--pedestal $%
\Upsilon (I_{j}^{(r)})$ given by vertices of an open oriented path from the
left to the right boundary of $I_{j}^{(r)}(\omega)$. If $I_{j}^{(r)}(\omega)$
is the leftmost (resp. the rightmost) interval, an $r$--pedestal $\Upsilon
(I_{j}^{(r)})$ is given by vertices of an open oriented path which starts
from some vertex $x \in [-L, -L + 2\lfloor l_r^{\alpha^{\prime}/
\alpha}\rfloor]$ and ends at the right boundary of $I_{j}^{(r)}(\omega)$
(resp. starts from the left boundary of $I_{j}^{(r)}(\omega)$ and ends at
some vertex $x \in [L- 2\lfloor l_r^{\alpha^{\prime}/ \alpha}\rfloor, L]$). (%
$l_1$ being large, we may assume that the length of an $(r-1)$--block is
always bounded above by $\lfloor l_r^{\alpha^{\prime}/ \alpha}\rfloor$,
according to (\ref{comprimento}) for $r$ replaced by $r-1$.)

\item the collection $D^{(r)}(\omega)=\{D_{j}^{(r)}(\omega) \}_{j} $, where $%
D_{j}^{(r)}(\omega)$ is the set of labels of the defected $(r-1)$--blocks
which are contained in $I_{j}^{(r)}(\omega)$.
\end{itemize}

For $\omega $ fixed, the sequence of pairs 
\begin{equation*}
J_{L}^{(k-1)}(\omega )=\{(I^{(1)}(\omega ),D^{(1)}(\omega )),\dots
,(I^{(k-1)}(\omega ),D^{(k-1)}(\omega ))\},
\end{equation*}%
is called $(k-1)$--\textit{itinerary}, and $(I^{(r)},D^{(r)})$, is called
the $r$--th step of the itinerary, for $1\leq r\leq k-1$. We shall now see
how to define the $k$--blocks and the continuation to a $k$--itinerary. When 
$k=M$ we will end up with only one or two intervals.

\noindent \textbf{Construction of} $k$--\textbf{blocks.} For any $\omega$
and for each $z\in [-L,L]$ we set $j_{z}^{k}=\min \{j:z\in I_{j}^{(k-1)}\}$, 
$\hat{\jmath}_{i}^{k}=j_{x_{i}^{(k)}}^{k}$, cf. notation at the beginning of
this section, $i=-l_M/l_k,\dots, l_M/l_{k}-1$, and define the intervals: 
\begin{equation*}
\widetilde{I}_{i}^{(k)}=\bigcup_{s=\hat{\jmath}_{i}^{k}+1}^{\hat{\jmath}%
_{i+1}^{k}}I_{s}^{(k-1)}=:[a^{(k)}_i,a^{(k)}_{i+1}]
\end{equation*}
as well as the following partition of $\Omega_L$:

\begin{eqnarray}
G(\widetilde{I}_{i}^{(k)}) &=&\displaystyle \bigcap\limits_{s=\hat{\jmath}%
_{i}^{k}+1}^{\hat{\jmath} _{i+1}^{k}}G(I_{s}^{(k-1)}),  \notag \\
H_s(\widetilde{I}_{i}^{(k)}) &=&\displaystyle B(I_{s}^{(k-1)})\cap %
\displaystyle\bigcap\limits_{u=\hat{\jmath}_i^k+1,u\neq s}^{\hat{\jmath}%
_{i+1}^{k}}G(I_{u}^{(k-1)}) ,  \notag \\
H(\widetilde{I}_{i}^{(k)}) &=&\displaystyle \bigcup\limits_{s=\hat{\jmath}%
_{i}^{k}+1}^{\hat{\jmath}_{i+1}^{k}} H_s(\widetilde{I}_{i}^{(k)}) ,  \notag
\\
B(\widetilde{I}_{i}^{(k)}) &=&\Omega_L\setminus\left( G( \widetilde{I}%
_{i}^{(k)})\cup H(\widetilde{I}_{i}^{(k)})\right).  \label{eventosk}
\end{eqnarray}

\noindent \textbf{Adjustment.} Given $\omega \in \Omega_L$, consider all $i$
for which $H_s( \widetilde{I}_{i}^{(k)})$ occurs for $s$ such that the
distance of the defected $(k-1)$--block $I_{s}^{(k-1)}\subset \widetilde{I}%
_{i}^{(k)}$ to the left endpoint $a^{(k)}_i$ (right endpoint $a^{(k)}_{i+1}$%
, resp.) is less than $\lfloor l_{k}^{\alpha^{\prime} / \alpha}\rfloor$. If
this set is non-empty take arbitrarily any such $\widetilde{I}_{i}^{(k)}$.

When the selected $\widetilde{I}_{i}^{(k)}$ is the leftmost (resp. the
rightmost) interval in $[-L,L]$, and the defect stays at distance less than $%
\lfloor l_{k}^{\alpha^{\prime} / \alpha}\rfloor$ from $-L$ (resp. $L$), we
set ${I}_{i}^{(k)} = \widetilde{I}_{i}^{(k)}$, and say that $G({I}%
_{i}^{(k)}) $ occurs (or that ${I}_{i}^{(k)}$ is in $G$ state for this $%
\omega$). Otherwise, we then check if $\widetilde{I} _{i-1}^{(k)}$ (resp.\ $%
\widetilde{I}_{i+1}^{(k)}$) has a defected block $I_{r}^{(k-1)}$ at distance
at most $3\lfloor l_{k}^{\alpha ^{\prime }/\alpha }\rfloor$ from $a^{(k)}_i $
(resp. from $a^{(k)}_{i+1}$), and

1) If yes, then we consider a new interval $I_{i-1}^{(k)} = \widetilde
I_{i-1}^{(k)}\cup \widetilde I_i^{(k)}$ (respectively $I_{i}^{(k)}
=\widetilde I_{i}^{(k)}\cup \widetilde I_{i+1}^{(k)}$) and say that $%
B(I_{i-1}^{(k)} )$ (resp. $B(I_{i}^{(k)} )$) occurs, or that the
corresponding interval is in state $B$;

\smallskip

2) If not, then we consider two new intervals: 
\begin{eqnarray}
I_{i-1}^{(k)}=\bigcup_{s=\hat{\jmath}_{i-1}^{k}+1}^{j_{a_{i}^{(k)}-l_{k}^{%
\alpha ^{\prime }/\alpha }}^{k}-1}I_{s}^{(k-1)}, &\qquad \quad
&I_{i}^{(k)}=\bigcup_{s=j_{a_{i}^{(k)}-l_{k}^{\alpha ^{\prime }/\alpha
}}^{k}}^{\hat{\jmath}_{i+1}^{k}}I_{s}^{(k-1)}  \label{carnaval4} \\
\Big(\text{respectively,}\qquad I_{i}^{(k)}=\bigcup_{s=\hat{\jmath}%
_{i}^{k}+1}^{j_{a_{i+1}^{(k)}+l_{k}^{\alpha ^{\prime }/\alpha
}}^{k}}I_{s}^{(k-1)} &\qquad \quad
&I_{i+1}^{(k)}=\bigcup_{s=j_{a_{i+1}^{(k)}+l_{k}^{\alpha ^{\prime }/\alpha
}}^{k}+1}^{\hat{\jmath}_{i+2}^{k}}I_{s}^{(k-1)}\Big ).  \label{carnaval5}
\end{eqnarray}%
In the situation of (\ref{carnaval4}) we say that $H(I_{i}^{(k)})$ occurs,
and say that $G(I_{i-1}^{(k)})$, $H(I_{i-1}^{(k)})$, $B(I_{i-1}^{(k)})$
occurs according to the occurrence of the corresponding $G(\widetilde{I}%
_{i-1}^{(k)})$, $H(\widetilde{I}_{i-1}^{(k)})$, $B(\widetilde{I}%
_{i-1}^{(k)}) $ (resp. in the situation of (\ref{carnaval5}) we say that $%
H(I_{i-1}^{(k)})$ occurs, and say that $G(I_{i+1}^{(k)})$, $H(I_{i+1}^{(k)})$%
, $B(I_{i+1}^{(k)})$ occurs according to the occurrence of $G(\widetilde{I}%
_{i+1}^{(k)})$, $H(\widetilde{I}_{i+1}^{(k)})$, $B(\widetilde{I}%
_{i+1}^{(k)}) $).

Finally we set $I_{i}^{(k)}=\widetilde{I}_{i}^{(k)} $ if $\widetilde{I}%
_{i}^{(k)}$ was not involved in the adjustment and say $G(I_{i}^{(k)})$, $%
H(I_{i}^{(k)})$, $B(I_{i}^{(k)})$ occurs if the corresponding $G(\widetilde{I%
}_{i}^{(k)})$, $H( \widetilde{I}_{i}^{(k)})$, $B(\widetilde{I}_{i}^{(k)})$
does occur.

To conclude this step, we re-numerate the intervals from left to right as $%
I_j^{(k)}\, j=1,\dots $. If after this step we are still left with intervals 
$I_{i}^{(k)}$ for which $H(I_{i}^{(k)})$ occurs and its defected interval $%
I_{s}^{(k-1)}$ stays within distance $\lfloor l_{k}^{\alpha ^{\prime
}/\alpha }\rfloor $ from one of the endpoints of $I_{i}^{(k)}$, then we
repeat the above procedure. After finitely many steps of this adjustment
procedure all $I_i^{(k)}$ for which $H(I_i^{(k)})$ occurs have their
defected $(k-1)$--block at distance larger than $\lfloor
l_k^{\alpha^{\prime} / \alpha}\rfloor$ from the boundary of $I_i^{(k)}$.

Once the adjustments are completed, the final intervals, always re-numerated
from left to right as $I_{j}^{(k)},\;j=1,\dots $, are called $k$--blocks. We
then consider the collection $D^{(k)}=\{D_j^{(k)}\}_{j}$ where $D_j^{(k)}$
gives the labels of the defected $(k-1)$--blocks contained in $I_{j}^{(k)}$.

We can always write $I_{j}^{(k)}=\displaystyle\bigcup%
\nolimits_{s_{0}(j)}^{s_{1}(j)}I_{s}^{(k-1)}$. It is easy to check that the
procedure is well defined (measurable) and the validity of the following
recursive estimate: 
\begin{equation}
l_{k}-(2\lfloor l_{k}^{\alpha ^{\prime }/\alpha }\rfloor
+6l_{k-1})<|I_{j}^{(k)}|\leq 3l_{k}+6l_{k-1}.  \label{vno}
\end{equation}

\medskip

\noindent $k$\textbf{-Pedestals.} Given the $k$--itinerary we shall
associate to each $k$--block $I_{j}^{(k)}(\omega)$ a state $G$ or $B$, and
the blocks in state $G$ will have a $k$--pedestal, to be defined below. When 
$|D_j^{(k)}(\omega)|\ge 2$, the block is said to be in state $B$, and it has
no $k$--pedestal.

\begin{itemize}
\item When $D(I_{j}^{(k)})=\emptyset$, all its sub-blocks $I_{s}^{(k-1)}$
are in state $G$. In this case we define $\Upsilon (I_{j}^{(k)})= %
\displaystyle\bigcup\nolimits_{s_{0}(j)}^{s_{1}(j)}\Upsilon (I_{s}^{(k-1)})$.%
\newline

\item If $D(I_{j}^{(k)})=\{r\}$ and $I_{j}^{(k)}$ is not the leftmost (resp.
rightmost) interval in $[-L,L]$, we define $\Upsilon _{\mathcal{L}%
}(I_{j}^{(k)})=\displaystyle\bigcup\nolimits_{s_{0}(j)}^{r-1}\Upsilon
(I_{s}^{(k-1)})$ and $\Upsilon _{\mathcal{R}}(I_{j}^{(k)})=\displaystyle%
\bigcup\nolimits_{r+1}^{s_{1}(j)}\Upsilon (I_{s}^{(k-1)}),$ called left and
right pedestals\footnote{%
From the occurrence of $G(I_{s}^{(k-1)})$ for all other $(k-1)$--blocks
within $I_{j}^{(k)}$, we know there exists an open oriented path connecting
the left boundary of $I_{j}^{(k)}$ to the right boundary of $I_{r-1}^{(k-1)}$
and an open oriented path connecting the left boundary of $I_{r+1}^{(k-1)}$
to the right boundary of $I_{j}^{(k)}$. These paths are obtained by
concatenation of the corresponding $\Upsilon (I_{s}^{(k-1)}),\;{s_{0}(j)}%
\leq s\leq {r-1}$ and $r+1\leq s\leq s_{1}(j)$, respectively.} of $%
I_{j}^{(k)}$, and check if there exists $x\in \Upsilon _{\mathcal{L}%
}(I_{j}^{(k)})$ and $y\in \Upsilon _{\mathcal{R}}(I_{j}^{(k)})$ with $%
y-x\leq \lfloor l_{k}^{\alpha ^{\prime }/\alpha }\rfloor $ such that $\omega
_{\{x,y\}}=1$:

-- If yes, we say that $I_{j}^{(k)}$ is in state $G$, and if the pair $(x,y)$
with $x\in \Upsilon _{ \mathcal{L}}(I_{j}^{(k)}),y\in \Upsilon_{\mathcal{R}%
}(I_{j}^{(k)}),\;y-x \leq \lfloor l_{k}^{\alpha ^{\prime }/\alpha }\rfloor $%
, and $\omega _{\{x,y\}}=1$ is not unique, we choose one in an arbitrary
way, and, once $(x,y)$ is chosen, denote $\mathcal{D} (I_{j}^{(k)}) =
[x+1,y-1]$, and define 
\begin{equation*}
\Upsilon (I_{j}^{(k)})=\left( \Upsilon _{\mathcal{L}}(I_{j}^{(k)})\cup
\Upsilon _{\mathcal{R}}(I_{j}^{(k)})\right) \setminus \mathcal{D}%
(I_{j}^{(k)}).
\end{equation*}
-- If such an open edge $\{x,y\}$ does not exist we say that $I_{j}^{(k)}$
is in $B$ state.\newline

\item If $|D(I_{j}^{(k)})|=1$ and $I_{j}^{(k)}$ is the leftmost (resp.
rightmost) interval in $[-L,L]$ whose unique defected $(k-1)$--block $%
I_{r}^{(k-1)}$ stays within distance $\lfloor l_k^{\alpha^{\prime}/\alpha}
\rfloor$ from $-L$ (resp. $L$), then we say that $I_{j}^{(k)}$ is in state $%
G $ and we define its $k$--pedestal as $\Upsilon (I_{j}^{(k)})= \displaystyle%
\bigcup\nolimits_{r+1}^{s_{1}(j)}\Upsilon (I_{s}^{(k-1)})$ (resp. $\Upsilon
(I_{j}^{(k)})= \displaystyle\bigcup\nolimits_{s_{0}(j)}^{r-1}\Upsilon
(I_{s}^{(k-1)})$). \newline

\item Finally, if $|D(I_{j}^{(k)})|=1$ and $I_{j}^{(k)}$ is the leftmost
(resp. rightmost) interval in $[-L,L]$, but its unique defected $(k-1)$%
--block $I_{r}^{(k-1)}$ does not stay within distance $\lfloor
l_k^{\alpha^{\prime}/\alpha} \rfloor$ from $-L$ (resp. $L$), then we use the
same procedure as if $I_{j}^{(k)}$ were not an extremal $k$--block.
\end{itemize}


This completes the $k$--th step, associating with each itinerary $J^{(k-1)}$
its continuation with a random sequence of $k$--blocks $I^{(k)}= \{I^{(k)}_j
\}_{j}$, re-numerated from left to right. Moreover, with each $k$--block we
associate one of the states $G$ or $B$. \medskip

\noindent \textbf{Structure of pedestals.} First we state a simple geometric
property of pedestals, which will be used in estimating the conditional
probability that a $k$--block $I_{j}^{(k)}$ is in state $G$, given that $%
|D_{j}^{(k)}|=1$. Our goal is to show that there exists a positive constant $%
C\equiv C(\alpha ,\alpha ^{\prime })$ such that if a $k$--block, $k\geq 1$, $%
I^{(k)}=[s,s^{\prime }]$ contains only one defected $(k-1)$--block, here
denoted by $[a,a^{\prime }]$, with corresponding left and right pedestals $%
\Upsilon _{\mathcal{L}}$ and $\Upsilon _{\mathcal{R}}$, spanning from $s$ to 
$a$ and from $a^{\prime }$ to $s^{\prime }$, respectively, then 
\begin{equation}
\Big|\Upsilon _{\mathcal{L}}\cap \lbrack a-\lfloor l_{k}^{\alpha ^{\prime
}/\alpha }\rfloor ,\,a]\Big|\geq Cl_{k}^{\alpha ^{\prime }/\alpha }\quad 
\text{and}\quad \Big|\Upsilon _{\mathcal{R}}\cap \lbrack a^{\prime
},\,a^{\prime }+\lfloor l_{k}^{\alpha ^{\prime }/\alpha }\rfloor ]\Big|\geq
Cl_{k}^{\alpha ^{\prime }/\alpha }.  \label{bagunca}
\end{equation}%
\begin{figure}[tbp]
\centering
\includegraphics[scale=1.2]{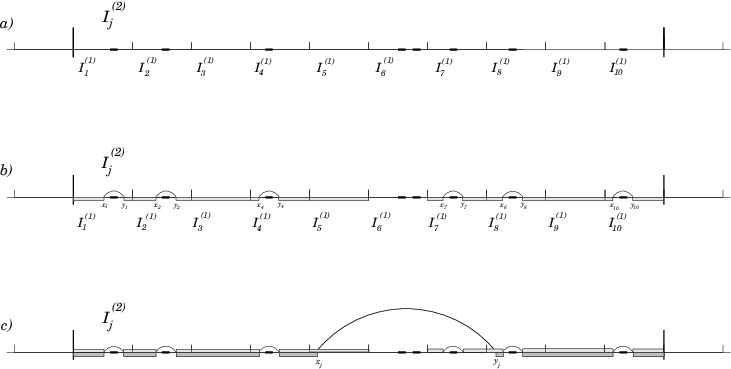}
\caption{Pedestals and defects: part $a)$ shows the deterministic $1$%
--blocks $I_{i}^{(1)},\;1\leq i\leq 10$, located in the $2$--block $%
I_{j}^{(2)}$; bold-face segments show location of the $0$--defects. Part $b)$
shows construction of $1$--pedestals, marked by light-gray strips. The block 
$I_{6}^{(1)}$ is a defected $1$--block. The segments $(x_{i},y_{i})$ are
\textquotedblleft enlarged\textquotedblright\ defects in $I_{i}^{(1)}$. Part 
$c)$ shows creation of $2$--pedestals, marked by dark-gray strips,
concatenated by long range edges. The segment $(x_{j},y_{j})$ is the
enlarged defect for $I_{j}^{(2)}$.}
\label{fig2}
\end{figure}
Inequality (\ref{bagunca}) follows trivially from the following recursive
relation: if we have a $k$--block $I^{(k)}=\cup _{s_{0}}^{s_{1}}I_{s}^{(k-1)}
$ which is in $G$ state, then 
\begin{equation*}
\big|\Upsilon (I^{(k)})\big|\geq \displaystyle\sum_{s\colon \lbrack
G(I_{s}^{(k-1)})\,\text{occurs}]}\big|\Upsilon (I_{s}^{(k-1)})\big|%
-l_{k}^{\alpha ^{\prime }/\alpha }.
\end{equation*}%
\smallskip 

We now give the announced basic estimate needed for the recursive step in
the previous construction, cf. (\ref{conserto}). Afterwards, we fix the
parameters which will determine the choice of $p$ close to one, as in (\ref%
{choice00}). In the lemma below, assume that $I^{(k)}_j$ is a $k$--block and 
$D^{(k)}_j=\{z\}$, i.e. the unique defected $(k-1)$--block within $I_j^{(k)}$
has index $z$, and by construction stays at distance larger than $\lfloor
l_{k}^{\alpha ^{\prime }/\alpha }\rfloor$ from the boundaries of $I_j^{(k)}$.

\begin{lemma}
\label{est33} There exists $\eta \equiv \eta (\alpha ,\alpha ^{\prime
},l_{1})$ with $\eta \searrow 0$ as $l_{1}\nearrow +\infty $ and such that
the following estimate for the conditional probability with respect to the
product measure (defined right above (\ref{occup})) 
\begin{equation}
\nu \big\lbrack \exists x\in \Upsilon _{\mathcal{L}}(I_{j}^{(k)}),y\in
\Upsilon _{\mathcal{R}}(I_{j}^{(k)}),\;y-x\leq \lfloor l_{k}^{\alpha
^{\prime }/\alpha }\rfloor :\omega _{\{x,y\}}=1\big||D^{(k)}_j|=1 \big\rbrack%
\geq 1-l_{k-1}^{-\beta (1-\eta )(\alpha ^{\prime }-1)}  \label{lax}
\end{equation}
holds for $k \ge 2$. For $k=1$ the r.h.s in (\ref{lax}) is replaced by $%
1-l_{1}^{-\beta (1-\eta )(\alpha ^{\prime }-1)/\alpha}$ .
\end{lemma}

\smallskip \noindent \textbf{Proof.} We show the above estimate by
conditioning on $D^{(k)}_j=\{z\}$, uniformly in $z$, and we make repeated
use of the following upper and lower bounds: if $I$ and $I^{\prime }$ are
two intervals, and $3\le d=\text{dist}(I,I^{\prime })$, then

\begin{equation}
C^{-}\ J(I,I^{\prime })\leq \sum_{\substack{ x\in I\cap \mathbb{Z}  \\ y\in
I^{\prime }\cap \mathbb{Z}}}\frac{1}{|x-y|^{2}}\leq C^{+}\ J(I,I^{\prime }),
\label{Bsumm}
\end{equation}%
holds with $C^{\pm }=\left( 1\pm 2/d\right) ^{2}$ and 
\begin{equation}
\ J(I,I^{\prime })=\int_{I\times I^{\prime }}dx\,dy\frac{1}{|x-y|^{2}}=\ln 
\frac{(|I|+d)(|I^{\prime }|+d)}{d\,(|I|+|I^{\prime }|+d)}.  \label{BJ}
\end{equation}%
Notice that we have $C^{-}\left( \left\vert x-y\right\vert -2\right)
^{-2}\leq \left\vert x-y\right\vert ^{-2}\leq C^{+}\left( \left\vert
x-y\right\vert +2\right) ^{-2}$ for $\left\vert x-y\right\vert \ge d$. We
shall need also the inequality%
\begin{equation}
J(I,I^{\prime })\leq 4\frac{\left\vert I^{\prime }\right\vert }{\left\vert
I^{\prime \prime }\right\vert }J\left( I,I^{\prime \prime }\right)
\label{inequal}
\end{equation}%
which holds for every $I$, $I^{\prime }$ and $I^{\prime \prime }$ such that $%
I^{\prime }\subset I^{\prime \prime }$ and $d^{\prime }=\text{dist}%
(I,I^{\prime \prime })\geq \left\vert I^{\prime \prime}\right\vert $.
Indeed, setting $f(x)=\int_{I}dy~|x-y|^{-2}$, for $x \in I^{\prime\prime}$,
straightforward calculations give that under the above conditions: 
\begin{equation*}
f(x^\prime)\le 4f(x^{\prime\prime}) \quad {\text{for each }}\quad x^\prime
\in I^{\prime },x^{\prime\prime} \in I^{\prime\prime}
\end{equation*}%
from where the inequality (\ref{inequal}) follows upon integration.

If $k \ge 2$ and $D^{(k)}_j=\{z\}$, we have the left $k$--pedestal $\Upsilon
_{\mathcal{L}}(I_{j}^{(k)})$ spanning from the left endpoint of $I_{j}^{(k)}$
to the left endpoint of $I_{z}^{(k-1)}$, and the right $k$ -pedestal $%
\Upsilon _{\mathcal{R}}(I_{j}^{(k)})$, spanning from the right endpoint of $%
I_{z}^{(k-1)}$ to the right endpoint of $I_{j}^{(k)}$. Take two segments $%
S_{z}^{\mathcal{L}}$ and $S_{z}^{\mathcal{R}}$, such that $|S_{z}^ {\mathcal{%
L}}|=|S_{z}^{\mathcal{R}}|=\lfloor\lfloor l_{k}^{\alpha ^{\prime }/\alpha
}\rfloor/3\rfloor$, lying immediately to the left and, respectively, to the
right of $I_{z}^{(k-1)}$. Denote 
\begin{eqnarray*}
\widehat{\Upsilon }_{\mathcal{L}}(I_{j}^{(k)}) &=&\Upsilon _{\mathcal{L}
}(I_{j}^{(k)})\cap S_{z}^{\mathcal{L}}, \\
\widehat{\Upsilon }_{\mathcal{R}}(I_{j}^{(k)}) &=&\Upsilon _{\mathcal{R}
}(I_{j}^{(k)})\cap S_{z}^{\mathcal{R}}.
\end{eqnarray*}
Then 
\begin{equation*}
\nu[{\text{all edges $\{ x,y\} ,\;x\in \widehat{\Upsilon }_{\mathcal{L}%
}(I_{j}^{(k)}),\,y\in \widehat{\Upsilon }_{\mathcal{R}}(I_{j}^{(k)})$ are
closed}}|\,D^{(k)}_j=\{z\}]
\end{equation*}%
\begin{equation}
\le{\prod_{\substack{ x\in S_{z}^{\mathcal{L}}  \\ y\in S_{z}^{\mathcal{R}}}}%
q_{\{ x,y\} }}{\prod_{\substack{ x\in S_{z}^{\mathcal{L}}\setminus \widehat{%
\Upsilon }_{\mathcal{L}}(I_{j}^{(k)})  \\ y\in S_{z}^{\mathcal{R}}}}q_{\{
x,y\} }^{-1}}{\prod _{\substack{ x\in S_{z}^{\mathcal{L}}  \\ y\in S_{z}^{%
\mathcal{R}}\setminus \widehat{\Upsilon }_{\mathcal{R}}(I_{j}^{(k)})}}q_{\{
x,y\} }^{-1}}.  \label{qqq}
\end{equation}

\noindent Applying (\ref{Bsumm}) to $S_{z}^{\mathcal{L}}$ and $S_{z}^{%
\mathcal{R}}$ we immediately get the following bound: 
\begin{equation}
{\prod_{\substack{ x\in S_{z}^{\mathcal{L}}  \\ y\in S_{z}^{\mathcal{R}}}}%
q_{\{ x,y\} }}=\exp \big\{-{\sum_{\substack{ x\in S_{z}^{\mathcal{L}}  \\ %
y\in S_{z}^{\mathcal{R}}}}\frac{\beta }{|x-y|^{2}}}\big\}\leq
l_{k-1}^{-\beta (\alpha ^{\prime }-1)(1-b)},  \label{gru}
\end{equation}
where $b\equiv b(\alpha ^{\prime },l_{1})$ and $b\searrow 0$ when $%
l_{1}\nearrow +\infty $. 
Similar computation gives that if a 1-block $I$ has a unique closed edge $\{
a, a+1\}$ with both $a, a+1$ at distance larger than $l_1^{\alpha^\prime/%
\alpha}$ from the endpoints of $I$, then the probability that there is an
open edge $\{ x,y\}$ with $x<a<y$, $y-x \le l_1^{\alpha^\prime/\alpha}$ is
larger than or equal of $1-l_1^{-\beta(1-\eta)(\alpha^\prime -1)/\alpha}$.

On the other hand denoting by $\mathcal{D}_{n}(S_{z}^{\mathcal{L}}\setminus 
\widehat{\Upsilon }_{\mathcal{L}}(I_{j}^{(k)})),\;0\leq n\leq k-2$ (resp. $%
\mathcal{D}_{n}(S_{z}^{\mathcal{R}}\setminus \widehat{\Upsilon }_{\mathcal{R}%
}(I_{j}^{(k)})))$ the set of vertices that belong to all defected $n$%
--blocks contained in the segment $S_{z}^{\mathcal{L}}$ (resp. $S_{z}^{%
\mathcal{R}}$), we get 
\begin{equation*}
{\prod_{\substack{ x\in S_{z}^{\mathcal{L}}\setminus \widehat{\Upsilon }_{%
\mathcal{L}}(I_{j}^{(k)}) \\ y\in S_{z}^{\mathcal{R}}}}q_{\{x,y\}}}%
=\prod_{n=0}^{k-2}{\prod_{\substack{ x\in \mathcal{D}_{n}(S_{z}^{\mathcal{L}%
}\setminus \widehat{\Upsilon }_{\mathcal{L}}(I_{j}^{(k)})) \\ y\in S_{z}^{%
\mathcal{R}}}}q_{\{x,y\}}}=\exp \Big\{-\sum_{n=0}^{k-2}\sum_{\substack{ x\in 
\mathcal{D}_{n}(S_{z}^{\mathcal{L}}\setminus \widehat{\Upsilon }_{\mathcal{L}%
}(I_{j}^{(k)})) \\ y\in S_{z}^{\mathcal{R}}}}\frac{\beta }{|x-y|^{2}}\Big\}.
\end{equation*}%
Once again, applying (\ref{Bsumm}) for each $0\leq n\leq k-2$ and taking
into account the structure of $n$--pedestals together with (\ref{inequal}),
we have (uniformly on all $l_{1}$ large enough) fixed positive constants $%
C_{i},i=1,2,3$ so that 
\begin{eqnarray*}
\sum_{n=0}^{k-2}\sum_{\substack{ x\in \mathcal{D}_{n}(S_{z}^{\mathcal{L}%
}\setminus \widehat{\Upsilon }_{\mathcal{L}}(I_{j}^{(k)})) \\ y\in S_{z}^{%
\mathcal{R}}}}\dfrac{\beta }{|x-y|^{2}} &\leq
&C_{1}\sum_{n=0}^{k-2}\sum_{\nu }J\left( I_{\nu }^{\prime },I^{\mathcal{R}%
}\right)  \\
&\leq &C_{2}\sum_{n=0}^{k-2}\frac{l_{n+1}^{\alpha ^{\prime }/\alpha }}{%
l_{n+1}}\sum_{\nu }J\left( I_{\nu }^{\prime \prime },I^{\mathcal{R}}\right) 
\\
&\leq &C_{3}l_{1}^{\alpha ^{\prime }/\alpha -1}J\left( I^{\mathcal{L}},I^{%
\mathcal{R}}\right) 
\end{eqnarray*}%
where $I_{\nu }^{\prime }$ and $I_{\nu }^{\prime \prime }$ are intervals in $%
\mathbb{R}$ so that $\bigcup_{\nu }\left( I_{\nu }^{\prime }\cap \mathbb{Z}%
\right) =\mathcal{D}_{n}(S_{z}^{\mathcal{L}}\setminus \widehat{\Upsilon }_{%
\mathcal{L}}(I_{j}^{(k)}))$, the sum $\sum_{\nu }$ is taken over all indices 
$\nu $ of $(n+1)$--blocks $I_{\nu }^{(n+1)}=:I_{\nu }^{\prime \prime }\cap 
\mathbb{Z}$ where the defected $n$--blocks are located, and moreover, $I^{%
\mathcal{L}}=\bigcup_{0\leq n\leq k-2}\bigcup_{\nu }I_{\nu }^{\prime \prime }
$ and  $I^{\mathcal{R}}$ is the convex envelop of $S_{z}^{\mathcal{R}}$. The
condition to apply (\ref{inequal}) in the first inequality above follows
from $3l_{k-1}+6l_{k-2}\leq l_{k}^{\alpha ^{\prime }/\alpha }$ which is true
for any $k\geq 2$, provided $l_{1}$ has been taken large enough. From this
we can easily get that 
\begin{equation}
{\prod_{\substack{ x\in S_{z}^{\mathcal{L}}\setminus \widehat{\Upsilon }_{%
\mathcal{L}}(I_{j}^{(k)}) \\ y\in S_{z}^{\mathcal{R}}}}q_{\{x,y\}}}\geq
l_{k-1}^{-\beta (\alpha ^{\prime }-1)b^{\prime }},  \label{hto}
\end{equation}%
where $b^{\prime }\equiv b^{\prime }(\alpha ,\alpha ^{\prime },l_{1})$ and $%
b^{\prime }\searrow 0$ when $l_{1}\nearrow +\infty $. Analogous lower bound
holds for the third term at the r.h.s of (\ref{qqq}). Finally, from the
upper bound for the length of a $(k-1)$-block, we have 
\begin{eqnarray}
&&[\omega :\exists x\in \Upsilon _{\mathcal{L}}(I_{j}^{(k)}),y\in \Upsilon _{%
\mathcal{R}}(I_{j}^{(k)}),\;y-x\leq \lfloor l_{k}^{\alpha ^{\prime }/\alpha
}\rfloor :\omega _{\{x,y\}}=1]^{c}  \notag \\
&\subseteq &[{\text{all edges $\{x,y\},\;x\in \widehat{\Upsilon }_{\mathcal{L%
}}(I_{j}^{(k)}),\,y\in \widehat{\Upsilon }_{\mathcal{R}}(I_{j}^{(k)})$ are
closed}}],~  \label{vno1}
\end{eqnarray}%
the statement of the Lemma follows from (\ref{gru}) and (\ref{hto}). $%
\square $

\medskip

\noindent \textbf{Fixing the parameters.} For fixed $\beta >1$, which is the
first main parameter of the model we choose the pair $\alpha ,\;\alpha
^{\prime }$ with $1<\alpha ^{\prime }<\alpha <2$ such that 
\begin{equation}
\beta (\alpha ^{\prime }-1)-\frac{2(\alpha -1)^{2}}{2-\alpha }>\alpha -1,
\label{aa}
\end{equation}%
i.e. $\beta (\alpha ^{\prime }-1)>\alpha \left( \alpha -1\right) /(2-\alpha
) $. We also fix 
\begin{equation}
\delta >\frac{2(\alpha -1)}{2-\alpha }  \label{del}
\end{equation}%
such that 
\begin{equation}
\beta (\alpha ^{\prime }-1)-\delta (\alpha -1)>\alpha -1.  \label{del2}
\end{equation}
By Lemma \ref{est33} we can fix $l_{1}>1$ so large that the parameter $\eta
=\eta (\alpha ,\alpha ^{\prime },l_{1})$ in (\ref{lax}) becomes so close to
zero, that 
\begin{equation}
\beta (1-\eta )(\alpha ^{\prime }-1)-\delta (\alpha -1)>\alpha -1.
\label{aald}
\end{equation}

\noindent Inequalities (\ref{lax}), (\ref{del}) and (\ref{aald}) are crucial
for the inductive estimates.

\medskip

\noindent \textbf{Cluster of the origin.} From the above estimates, and
recalling (\ref{vno}), the initial heuristic discussion is indeed made
rigorous: for the above choice of parameters and picking $l_{1}$ large
enough we (recursively) obtain that for all $M\ge 1$ and at all scales $%
k=1,\dots ,M$, 
\begin{equation}
\nu (I_{j}^{(k)}\text{ is defected })\leq l_{k}^{-\delta }.  \label{psiu3}
\end{equation}%
Indeed, due to (\ref{vno}), we see that the previous analysis and the above
choice of the parameters turns rigorous the discussion leading to (\ref{psiu}%
) and (\ref{aaa}). Now, for $k=M$, we have at most two $M$--blocks, denoted
by $I_{i}^{(M)}$, where $1\leq i\leq s$ and $s(\omega )\in \{1,2\}$. In
particular, from (\ref{psiu3}), we immediately have the basic estimate (\ref%
{basica}) announced in the introduction. Next we give the uniform lower
bound for 
\begin{equation*}
\nu \left( 0\rightsquigarrow y,\;{\text{for some}}\;y\in \big[L-2\lfloor l_{{%
M}}^{\alpha ^{\prime }/\alpha }\rfloor ,L\big]\right) .
\end{equation*}

Recalling that $j_{z}^{k}=\min \{j:z\in I_{j}^{(k-1)}\}$, for any $1\leq
k\leq {M}$ we define the following events:

\begin{equation}
\psi ^{(k)}=\bigcap_{i=j_{0}^{k}-\lfloor \lfloor l_{k}^{\alpha ^{\prime
}/\alpha }\rfloor /l_{k-1}\rfloor }^{j_{0}^{k}+\lfloor \lfloor l_{k}^{\alpha
^{\prime }/\alpha }\rfloor /l_{k-1}\rfloor }G(I_{i}^{(k-1)})
\label{eventos1}
\end{equation}
and consider 
\begin{equation}
\Psi_{M}=\bigcap_{j=1}^{{M}}\psi ^{(j)}.  \label{eventos2}
\end{equation}

The occurrence of $\bigcap_{k=1}^{n}\psi ^{(k)}$, $1\leq n\leq M$ implies
that the origin $0$ is the right (resp. left) end--vertex of a $(n-1)$%
--block $I_{j_{0}^{n}}^{(n-1)}$ (resp. $I_{j_{0}^{n}+1}^{(n-1)}$) for each $%
n $, since no adjustments are performed in this case, and necessarily it
belongs to the pedestals $\Upsilon (I_{j_{0}^{n}}^{(n-1)})$ and $\Upsilon
(I_{j_{0}^{n}+1}^{(n-1)})$ for any $\;1\leq n\leq {M}$. In particular, for $%
\omega \in \Psi_{M}$ we have $s(\omega )=2$. Moreover, in the event $\Psi
_{M}\cap G(I_{1}^{(M)})\cap G(I_{2}^{(M)})$, the origin 0 belongs to an open
oriented path connecting $\big[-L, -L+2\lfloor l_{{M}}^{\alpha ^{\prime
}/\alpha }\rfloor\big]$ to $\big[L-2\lfloor l_{{M}}^{\alpha ^{\prime
}/\alpha }\rfloor ,L\big]$ as described above.

\noindent Taking into account the estimate (\ref{psiu3}) and the definition (%
\ref{eventos1}) we have for $k \ge 2$: 
\begin{equation*}
\nu(\psi^{(k)}) \ge 1-
(2(l_{k-1})^{\alpha^{\prime}-1}+1)(l_{k-1})^{-\delta}\ge 1-
3(l_{k-1})^{\alpha^{\prime}-1}(l_{k-1})^{-\delta}.
\end{equation*}
Since 
\begin{equation*}
\delta -(\alpha ^{\prime }-1)>\delta -(\alpha -1)\geq \frac{\alpha (\alpha
-1)}{2-\alpha }>0,
\end{equation*}
we define $u= \delta - (\alpha^{\prime}-1) >0$ and rewrite the above
inequality: 
\begin{equation*}
\nu(\psi^{(k)}) \ge 1- 3(l_{k-1})^{-u} \qquad \text{for} \; k \ge 2.
\end{equation*}

\noindent Since $l_k$ grow super-exponentially fast, we get immediately that
the series 
\begin{equation*}
(l_1)^{-u} +(l_2)^{-u}+ (l_3)^{-u} +... = S (l_1)
\end{equation*}
converges and 
\begin{equation*}
S(l_1)\longrightarrow 0, \; \text{when} \; l_1 \to \infty.
\end{equation*}
This immediately implies that 
\begin{equation}
\nu \big( \big[ \bigcap_{j=2}^{{M}}\psi ^{(j)} \cap G(I_1^{M}) \cap G(I_2^M)%
\big]^c \big)  \label{eventos2a}
\end{equation}
can be made arbitrarily small, uniformly in $M$.

\noindent Finally, by choosing $l_1$ large enough, and then $p$ close enough
to $1$ we get that $\nu (\psi^{(1)})$ can be made arbitrarily close to $1$.


\noindent The proof of Theorem \ref{percolation} follows at once. $\Box $




\bigskip

\noindent \textbf{Acknowledgments} The authors thank C. M. Newman for
suggesting this problem and A.-S. Sznitman for many useful discussions. We
thank E. Presutti for pointing out the \textquotedblleft interface
question\textquotedblright\ and M. Cassandro and I. Merola for discussions
on this topic. We thank the referees of a previous version, for their
careful reading and for their suggestions. This work was supported by Faperj
grant E-26/151.905/2000, Fapesp grant 03/01366-7, CNPq grant 477259/01-4 and
Fundacion Andes. VS is partially supported by CNPq grant 302221/2008-5. MEV
is partially supported by CNPq grant 302796/2002-9. We also thank various
institutions CBPF, IMPA, IF-USP, ETH, UC-Berkeley and MSRI for financial
support and hospitality during periods when this work was carried out.

\end{document}